\documentclass[12pt, twoside, leqno]{article}
\usepackage{amsmath,amsthm}
\usepackage{amssymb,latexsym}
\usepackage{enumerate}

\newtheorem{thm}{Theorem}[section]
\newtheorem{prop}[thm]{Proposition}
\newtheorem{cor}[thm]{Corollary}
\newtheorem{lem}[thm]{Lemma}
\newtheorem{prob}[thm]{Problem}

\theoremstyle{definition}

\newtheorem{rem}[thm]{Remark}

\numberwithin{equation}{section}

\usepackage{mathrsfs}
\usepackage{CJK}
\usepackage{amsmath, amscd}
\usepackage{fancyhdr}
\usepackage{enumerate}

\begin{document}

\baselineskip=17pt

\title{The kernels of the linear maps of finite Abelian group algebras}
\author{ Dan Yan \footnote{ The author is supported by the NSF of China (Grant No. 11871241; 11601146), the China Scholarship Council and the Construct Program of the Key Discipline in Hunan Province.}\\
MOE-LCSM,\\ School of Mathematics and Statistics,\\
 Hunan Normal University, Changsha 410081, China \\
\emph{E-mail:} yan-dan-hi@163.com \\
}
\date{}

\maketitle

\renewcommand{\thefootnote}{}

\renewcommand{\thefootnote}{\arabic{footnote}}
\setcounter{footnote}{0}

\begin{abstract} In our paper, we give a necessary and sufficient conditions for the kernels of the linear maps of finite Abelian group algebras to be Mathieu-Zhao spaces of $K[G]$ if $G$ is a finite Abelian group and $K$ is a split field for $G$. Hence we classify all Mathieu-Zhao spaces of the finite Abelian group algebras if $K$ is a split field for $G$.
\end{abstract}
 \vskip 2.5mm

\section{Introduction}

Throughout this paper, we will write $K$ for a field without specific note and $K[G]$ for the group algebra of $G$ over $K$. $V_G$ is the $K$-subspace of the group algebra $K[G]$ consisting of all the elements of $K[G]$ whose coefficient of the identity element $1_G$ of $G$ is equal to zero. $L$ is a linear map from $K[G]$ to $K$ and $L|_H$ means restricting $L$ on $H$, where $H$ is a subgroup of $G$. $H$ is a $p'$-subgroup of $G$ if $p\nmid |H|$.

The Mathieu-Zhao space was introduced by Wenhua Zhao in \cite{12}, which is a natural generalization of ideals. The notion was named after A. van den Essen's suggestion and Mathieu conjecture. In \cite{5}, J. J. Duistermaat and W. van der Kallen proved the Mathieu conjecture for the case of tori, which can be re-stated as follows.
\begin{thm}
Let $z=(z_1,z_2,\ldots,z_m)$ be $m$ commutative free variables and $V$ the subspace of the Laurent polynomial algebra $\mathbb{C}[z^{-1},z]$ consisting of the Laurent polynomials with no constant term. Then $V$ is a Mathieu-Zhao space of $\mathbb{C}[z^{-1},z]$.
\end{thm}
Let $G$ be the free Abelian group $\mathbb{Z}^m$ $(m\geq 1)$. Then the Laurent polynomial algebra $\mathbb{C}[z^{-1},z]$ can be identified with the group algebra $\mathbb{C}[G]$. Under this identification, the subspace of $V$ in the theorem is $V_G$. In \cite{2}, Wenhua Zhao and R. Willems proved that $V_G$ is a Mathieu-Zhao space of $K[G]$ if $G$ is a finite group and $\operatorname{char}K=0$ or $\operatorname{char}K=p>|G|$. For finite Abelian group, they proved that if $K$ contains a primitive $d$-th root of unity, then $V_G$ is a Mathieu-Zhao space of $K[G]$ if and only if $\operatorname{char}K=p>d$, where $|G|=p^ad$, $p\nmid d$. In \cite{19}, Wenhua Zhao and the author give a sufficient and necessary condition for $V_G$ to be a Mathieu-Zhao space of $K[G]$ if $G$ is a finite group and $K$ is a split field for $G$. Hence it's natural to ask the following question.
\begin{prob}\label{prob1.2}
Let $G$ be a finite group with $|G|=n$ and $L=(L_1,L_2,\ldots,L_r)$, $L_i$ is a linear map from $K[G]$ to $K$ such that $L_i(g_j)=l_{i,j}$ for all $1\leq i\leq r$, $1\leq j\leq n$. Suppose that $L_1,L_2,\ldots,L_r$ are linearly independent over $K$. Then under what conditions on $L$ and $K$, $\operatorname{Ker}L$ forms a Mathieu-Zhao space of the group algebra $K[G]$?
\end{prob}
It's easy to see that if $r\geq n$, then $\operatorname{Ker}L=0$. If $r\leq n-1$, then $\operatorname{dim}_K\operatorname{Ker}L=n-r$ and every codimension $r$ subspace of $K[G]$ is $\operatorname{Ker}L$ for some linear map $L$. Hence $\operatorname{Ker}L$ are all the codimension $r$ subspaces of $K[G]$.

In our paper, we first prove some properties of $\operatorname{Ker}L$ in section 2. Then we classify all Mathieu-Zhao spaces of $K[G]$ if $G$ is a finite Abelian group and $K$ is a split field for $G$ in section 3. Thus, we solve Problem \ref{prob1.2} if $G$ is a finite Abelian group.

\section{Some properties of Ker$L$}

\begin{prop}\label{prop2.1}
Let $K$, $G$, $L$ be as in Problem \ref{prob1.2} and $g_1$ is the identity $1_G$ of $G$. Then we have the following statements:

$(1)$ If all the $l_{i,j}$ are equal for all $1\leq i\leq r$, $1\leq j\leq n$, then $\operatorname{Ker}L$ is an ideal of $K[G]$.

$(2)$ If $\operatorname{Ker}L$ is a Mathieu-Zhao space of $K[G]$, then there exists $i_0\in\{1,2,\ldots,r\}$ such that $l_{i_0,1}\neq 0$.
\end{prop}
\begin{proof}
 $(1)$ Let $l:=l_{i,j}$ for all $1\leq i\leq r$, $1\leq j\leq n$. Then $\operatorname{Ker}L=\{\alpha=\sum_{j=1}^nc_jg_j\in K[G]|l\cdot\sum_{j=1}^nc_j=0\}$. Since $l\neq 0$, we have that $\operatorname{Ker}L=\{\alpha=\sum_{j=1}^nc_jg_j\in K[G]|\sum_{j=1}^nc_j=0\}$. It's easy to check that $\operatorname{Ker}L$ is an ideal of $K[G]$.

$(2)$ If $l_{1,1}=\cdots=l_{r,1}=0$, then $1_G\in \operatorname{Ker}L$. If $\operatorname{Ker}L$ is a Mathieu-Zhao space of $K[G]$, then $\operatorname{Ker}L=K[G]$. That is, $L=0$, which is a contradiction. Then the conclusion follows.
\end{proof}

\begin{rem}\label{rem2.2}
We can see from Proposition \ref{prop2.1} that we can assume $l_{i_0,1}\neq 0$ for some $i_0\in \{1,2,\ldots,r\}$ in the following arguments. If $r=1$ and $l_{1,2}=l_{1,3}=\cdots=l_{1,n}=0$, $l_{1,1}\neq 0$, then $\operatorname{Ker}L=V_G$, which is discussed in \cite{2} and \cite{19}.
\end{rem}

\begin{prop}\label{prop2.3}
Let $R$ be any commutative ring and $G$ any group. Suppose that $L=(L_1,L_2,\ldots,L_r)$ is a linear map from $R[G]$ to $R$. If $\operatorname{Ker}L$ is a Mathieu-Zhao space of $R[G]$, then $\operatorname{Ker}(L|_H)$ is a Mathieu-Zhao space of $R[H]$, where $H$ is any subgroup of $G$.
\end{prop}
\begin{proof}
Assume otherwise. Then there exist $u,v_1,v_2\in R[H]$ such that $u^m\in \operatorname{Ker}(L|_H)$ for all $m\geq 1$, but $v_1u^mv_2\notin \operatorname{Ker}(L|_H)$ for infinitely many $m\geq 1$.

Since $R[H]\subseteq R[G]$, we have $u,v_1,v_2\in R[G]$ and $u^m\in \operatorname{Ker}L$ for all $m\geq 1$, but $v_1u^mv_2\notin \operatorname{Ker}L$ for infinitely many $m\geq 1$. Otherwise, $v_1u^mv_2\in \operatorname{Ker}L\cap R[H]=\operatorname{Ker}(L|_H)$, which is a contradiction. Hence $\operatorname{Ker}L$ is not a Mathieu-Zhao space of $R[G]$, which is a contradiction. Then the conclusion follows.
\end{proof}

\begin{cor}\label{cor2.4}
Let $L$, $G$ be as in Problem \ref{prob1.2} and $K$ be a field of characteristic $p$, $H$ a normal subgroup of $G$. If $H$ is a $p'$-subgroup and $\operatorname{Ker}L$ is a Mathieu-Zhao space of $K[G]$, then $\operatorname{Ker}(L|_{G/H})$ is a Mathieu-Zhao space of $K[G/H]$.
\end{cor}
\begin{proof}
Let $\varphi$ be the natural surjective homomorphism from $K[G]$ to $K[G/H]$ and $E_H=\frac{1}{|H|}\sum_{j=1}^{|H|}h_j$. Then $(1-E_H)K[G]=\operatorname{Ker}\varphi$ and $E_HK[G]\cong K[G/H]$. Thus, we have $K[G]\cong(1-E_H)K[G]\oplus K[G/H]$. Therefore, $K[G/H]$ can be seen as a subalgebra of $K[G]$. It follows from the arguments of Proposition \ref{prop2.3} that $\operatorname{Ker}(L|_{G/H})$ is a Mathieu-Zhao space of $K[G/H]$.
\end{proof}

\begin{lem}\label{lem2.5}
Let $L$ and $G$ be as in Problem \ref{prob1.2}. Then $\operatorname{Ker}L=\{\beta\in K[G]|\operatorname{Tr}\beta\alpha_i\allowbreak=0~for~all~1\leq i\leq r\}$, where $\alpha_i=\sum_{j=1}^nl_{i,j}g_j^{-1}$ for $1\leq i\leq r$.
\end{lem}
\begin{proof}
Let $\beta=\sum_{j=1}^nc_jg_j$. Then $L_i(\beta)=\sum_{j=1}^nc_jl_{i,j}=\operatorname{Tr}\beta\alpha_i$ for all $1\leq i\leq r$. Hence the conclusion follows.
\end{proof}

\begin{thm}\label{thm2.6}
Let $L$ and $G$ be as in Problem \ref{prob1.2} and $K$ a field of characteristic zero or a field of characteristic $p$ and $p\nmid |G|$. If $K$ is a split field for $G$, then
$$\operatorname{Ker}L\cong\{(A_1,\ldots,A_s)\in A|\sum_{j=1}^sn_j\operatorname{Tr}(C_{i,j}A_j)=0~for~all~1\leq i\leq r\},$$
where $A=M_{n_1}(K)\times \cdots \times M_{n_s}(K)$ is the product of matrices and $C_{i,j}=\rho_j(\alpha_i)\in M_{n_j}(K)$ for $1\leq j\leq s$ and for all $1\leq i\leq r$, where $\alpha_i$ be as in Lemma \ref{lem2.5} for $1\leq i\leq r$, $\rho_j$ is an irreducible representation of $G$ and $n_j=\rho_j(1)$ for $1\leq j\leq s$ and $s$ is the number of distinct (up to isomorphism) irreducible representations of $G$.
\end{thm}
\begin{proof}
Since $\operatorname{char}K=0$ or $\operatorname{char}K=p$ and $p\nmid |G|$, we have that $K[G]$ is semi-simple. Since $K$ is a split field for $G$, we have that
$$K[G]\cong M_{n_1}(K)\times M_{n_2}(K)\times \cdots\times M_{n_s}(K),$$
where $M_{n_j}(K)$ is the ring of $n_j\times n_j$ matrices over $K$ for $1\leq j\leq s$.
Let $\tilde{\rho}$ be the regular representation of $K[G]$. Then $\operatorname{Tr}(\beta)=0$ if and only if $\operatorname{Tr}(\tilde{\rho}(\beta))=0$ for all $\beta\in K[G]$. Let $\rho=(\rho_1,\rho_2,\ldots,\rho_s)$. Then $\rho$ is a ring isomorphism from $K[G]$ to $A$. Let $\beta$ be any element in $K[G]$. Then
$$\rho(\alpha_i\beta)=(\rho_1(\alpha_i\beta),\rho_2(\alpha_i\beta),\ldots,\rho_s(\alpha_i\beta))=(\rho_1(\alpha_i)\rho_1(\beta),\ldots,\rho_s(\alpha_i)\rho_s(\beta)).$$
Suppose that
$$\rho(\alpha_i)=(\rho_1(\alpha_i),\ldots,\rho_s(\alpha_i))=\left( \begin{matrix}
C_{i,1} & 0 &\cdots & 0 \\
0 & C_{i,2} & \cdots & 0 \\
\vdots & \vdots &\ddots& \vdots\\
0 & 0 &\cdots & C_{i,s}
\end{matrix} \right)\in A$$
and
$$\rho(\beta)=(\rho_1(\beta),\ldots,\rho_s(\beta))=\left( \begin{matrix}
A_1 & 0 &\cdots & 0 \\
0 & A_2 & \cdots & 0 \\
\vdots & \vdots &\ddots& \vdots\\
0 & 0 &\cdots & A_s
\end{matrix} \right)\in A$$
for all $1\leq i\leq r$.
Then we have that
$$\rho(\alpha_i\beta)=\left( \begin{matrix}
C_{i,1}A_1 & 0 &\cdots & 0 \\
0 & C_{i,2}A_2 & \cdots & 0 \\
\vdots & \vdots &\ddots& \vdots\\
0 & 0 &\cdots & C_{i,s}A_s
\end{matrix} \right)\in A.$$
Thus, we have the following commutative diagram:
\[\begin{CD}
  K[G] @  >\cong>> M_{n_1}(K)\times M_{n_2}(K)\times \cdots\times M_{n_s}(K) \\
  @V \tilde{\rho}(\alpha_i\beta) VV @V \phi(\rho(\alpha_i\beta)) VV  \\
  K[G] @>\cong>> \left( \begin{matrix}
M_{n_1}(K) & 0 &\cdots & 0 \\
0 & M_{n_2}(K) & \cdots & 0 \\
\vdots & \vdots &\ddots& \vdots\\
0 & 0 &\cdots & M_{n_s}(K)\end{matrix} \right),
\end{CD}
\]
where $\phi$ is the natural isomorphism between the two algebras. Thus, we have
that $\operatorname{Tr}(\tilde{\rho}(\alpha_i\beta))=0$ if and only if $\operatorname{Tr}(\phi(\rho(\alpha_i\beta)))=0$. Since $\operatorname{Tr}(\phi(\rho(\alpha_i\beta)))=n_1\operatorname{Tr}(C_{i,1}A_1)+n_2\operatorname{Tr}(C_{i,2}A_2)+\cdots+n_s\operatorname{Tr}(C_{i,s}A_s)$, we have that $\operatorname{Tr}(\alpha_i\beta)=0$ if and only if $n_1\operatorname{Tr}(C_{i,1}A_1)+n_2\operatorname{Tr}(C_{i,2}A_2)+\cdots+n_s\operatorname{Tr}(C_{i,s}A_s)=0$ for all $1\leq i\leq r$.
Thus, we have that $\operatorname{Ker}L\cong V$, where
$$V=\{(A_1,A_2,\ldots,A_s)\in A|\sum_{j=1}^sn_j\operatorname{Tr}C_{i,j}A_j=0~\operatorname{for}~\operatorname{all}~1\leq i\leq r\}.$$
\end{proof}

\begin{cor}\label{cor2.7}
Let $L$ and $G$ be as in Problem \ref{prob1.2} and $K$ a field of characteristic zero or a field of characteristic $p$ and $p\nmid |G|$. If $K$ is a split field for $G$ and $r=1$, then $\operatorname{Ker}L$ is a Mathieu-Zhao space of $K[G]$ if and only if $n_1\lambda_1d_1+n_2\lambda_2d_2+\cdots+n_t\lambda_td_t\neq 0$ for all non-zero vectors $\tilde{d}=(d_1,\ldots,d_t)$, $d_j\in \{0,1,\ldots,n_j\}$ for $1\leq j\leq t$, where $n_j\lambda_j=\operatorname{Tr}\rho_j(\alpha_1)$ and $\alpha_1$ be as in Lemma \ref{lem2.5}, $\rho_j$ is an irreducible representation of $G$ for $1\leq j\leq s$ and $s$ is the number of distinct (up to isomorphism) irreducible representations of $G$ and $t\in\{1,2,\ldots,s\}$.
\end{cor}
\begin{proof}
It follows from Theorem \ref{thm2.6} that $\operatorname{Ker}L\cong V$, where $V=\{(A_1,\ldots,A_s)\in A|\sum_{j=1}^sn_j\operatorname{Tr}(C_{1,j}A_j)=0\}$ and $C_{1,j}=\rho_j(\alpha_1)\in M_{n_j}(K)$ for $1\leq j\leq s$. Let $\rho$ be as in Theorem \ref{thm2.6}. Since $\alpha_1\neq 0$ and $\rho$ is an isomorphism, we have that $\rho(\alpha_1)\neq 0$. We can assume that $C_{1,1},\ldots,C_{1,t}$ are not equal to zero and $C_{1,t+1}=\cdots=C_{1,s}=0$ for some $t\in \{1,2,\ldots,s\}$ by reordering the $\rho_j$ for $1\leq j\leq s$. It follows from Theorem 5.8.1 in \cite{3} or Theorem 4.4 in \cite{6} that $V$ is a Mathieu-Zhao space of $A$ if and only if $C_{1,j}=\lambda_jI_{n_j}$ and $n_1\lambda_1d_1+\cdots+n_t\lambda_td_t\neq 0$ for all nonzero vectors $\tilde{d}=(d_1,\ldots,d_t)$ and $d_j\in \{0,1,\ldots,n_j\}$ for $1\leq j\leq t$. Then the conclusion follows.
\end{proof}

\begin{cor}\label{cor2.8}
Let $L$ and $G$ be as in Problem \ref{prob1.2} and $K$ a field of characteristic zero or a field of characteristic $p$ and $p\nmid |G|$. If $K$ is a split field for $G$ and $r=1$, then the following two statements are equivalent:

$(1)$ $\operatorname{Ker}L$ is a Mathieu-Zhao space of $K[G]$.

$(2)$ There exists $\mu_1,\ldots,\mu_t\in K$ such that $L_1=\mu_1\chi_1+\mu_2\chi_2+\cdots+\mu_t\chi_t$ and $\mu_1d_1+\cdots+\mu_td_t\neq 0$ for all nonzero vectors $\tilde{d}=(d_1,d_2,\ldots,d_t)$, $d_j\in\{0,1,\ldots,n_j\}$ for $1\leq j\leq t$, where $\chi_1, \chi_2,\ldots,\chi_s$ are the non-isomorphic irreducible characters of $G$ and $\mu_j=n^{-1}n_j\lambda_j$, $n_j=\chi_j(1)$, $n_j\lambda_j=\operatorname{Tr}\rho_j(\alpha_1)$, $\alpha_1$ be as in Lemma \ref{lem2.5} and $\rho_j$ is an irreducible representation of $G$ with character $\chi_j$ for $1\leq j\leq s$, $s$ is the number of distinct (up to isomorphism) irreducible representations of $G$ and $t\in\{1,2,\ldots,s\}$.
\end{cor}
\begin{proof}
$(1)\Rightarrow (2)$
Since $L_1(\beta)=\operatorname{Tr}(\alpha_1\beta)$ for any $\beta\in K[G]$, where $\alpha_1$ be as in Lemma \ref{lem2.5}, we have that
$$n\operatorname{Tr}(\alpha_1\beta)=\operatorname{Tr}\tilde{\rho}(\alpha_1\beta)=\operatorname{Tr}\phi (\rho(\alpha_1\beta))$$
by following the arguments of Theorem \ref{thm2.6}, where $\tilde{\rho}$ be as in Theorem \ref{thm2.6}. Since $\operatorname{Ker}L$ is a Mathieu-Zhao space of $K[G]$, it follows from Corollary \ref{cor2.7} that $C_{1,j}=\lambda_jI_{n_j}$ for $\lambda_j\in K$ and for all $1\leq j\leq s$. We can assume that $\lambda_1\cdots \lambda_t\neq 0$ and $\lambda_{t+1}=\cdots=\lambda_s=0$ for some $t\in \{1,2,\ldots,s\}$ by reordering $\chi_1, \chi_2, \ldots, \chi_s$.

Thus, it follows from Lemma \ref{lem2.5} that
$$L_1(\beta)=\operatorname{Tr}(\alpha_1\beta)=n^{-1}(n_1\lambda_1\operatorname{Tr}A_1+n_2\lambda_2\operatorname{Tr}A_2+\cdots+n_t\lambda_t\operatorname{Tr}A_t).$$
Since $\operatorname{Tr}A_j=\chi_j(\beta)$ for all $1\leq j\leq s$, we have that
$$L_1=n^{-1}(n_1\lambda_1\chi_1+n_2\lambda_2\chi_2+\cdots+n_t\lambda_t\chi_t).$$
It follows from Corollary \ref{cor2.7} that $n_1\lambda_1d_1+\cdots+n_t\lambda_td_t\neq 0$ for all nonzero vectors $\tilde{d}=(d_1,d_2,\ldots,d_t)$, $d_j\in\{0,1,\ldots,n_j\}$ for $1\leq j\leq t$. Let $\mu_j=n^{-1}n_j\lambda_j$ for all $1\leq j\leq s$. Then the conclusion follows.

$(2)\Rightarrow (1)$ Since $\operatorname{Ker}L=\{\beta\in K[G]|L_1(\beta)=0\}=\{\beta\in K[G]|\mu_1\chi_1(\beta)+\cdots+\mu_t\chi_t(\beta)=0\}$ and there exists $A_j\in M_{n_j}(K)$ such that $\operatorname{Tr}A_j=\chi_j(\beta)$ for all $1\leq j\leq t$, we have that
$$\operatorname{Ker}L=\{(A_1,\ldots,A_t)\in M_{n_1}(K)\times\cdots\times M_{n_t}(K)|\sum_{j=1}^t\mu_j\operatorname{Tr}A_j=0\}.$$ Then the conclusion follows from Theorem 5.8.1 in \cite{3} or Theorem 4.4 in \cite{6}.
\end{proof}

\begin{rem}\label{rem2.9}
To prove that if $n_1\lambda_1d_1+n_2\lambda_2d_2+\cdots+n_t\lambda_td_t\neq 0$ for all nonzero vectors $\tilde{d}=(d_1,d_2,\ldots,d_t)$ and $d_j\in \{0,1,\ldots,n_j\}$ for $1\leq j\leq t$, then $\operatorname{Ker}L$ is a Mathieu-Zhao space of $K[G]$ for $r=1$, we don't need the condition that $K$ is a split field for $G$ in Corollary \ref{cor2.7} by following the arguments Theorem 5.8.1 in \cite{3} because an idempotent matrix can be conjugated to a diagonal matrix with only 0 and 1 on the diagonal over division rings.

If $L=\mu_j\chi_j$ for some $j\in \{1,2,\ldots,t\}$, $\mu_j\in K^*$, then it follows from the arguments of Corollary \ref{cor2.8} that the condition $n_1\lambda_1d_1+n_2\lambda_2d_2+\cdots+n_t\lambda_td_t \neq 0$ in Theorem \ref{thm2.6}
is equivalent to $n_jd_j\neq 0$ for all $1\leq d_j\leq n_j$, which is clearly true if $\operatorname{char}K=0$. If $\operatorname{char}K=p$, then the condition is equivalent to $p>n_j$. To see this, we can assume that $p|n_jd_j$ for some $d_j\in \{1,2,\ldots,n_j\}$, then $p|n_j$ or $p|d_j$, which contradicts with $p>n_j$. Thus, if $p>n_j$, then $n_jd_j\neq 0$ mod $p$ for all $1\leq d_j\leq n_j$. Conversely, suppose that $p\leq n_j$. Then let $d_j=p\in \{1,2,\ldots,n_j\}$, we have that $n_jp=0$ mod $p$, which is a contradiction. Hence if $n_jd_j\neq 0$ mod $p$ for all $1\leq d_j\leq n_j$, then $p>n_j$. Therefore, the conclusion is the same as Theorem 5.1 in \cite{17} in this situation.
\end{rem}

\section{Mathieu-Zhao spaces of finite Abelian group algebras}

\begin{prop}\label{prop3.1}
Let $B=K\times\cdots\times K$ be a $K$-algebra and
$$V=\{(a_1,a_2,\ldots,a_n)\in B|\sum_{j=1}^n\gamma_{i,j}a_j=0~for~all~1\leq i\leq r\},$$
where $\gamma_{i,j}\in K$ for all $1\leq i\leq r$, $1\leq j\leq n$.
If at least one of $\gamma_{i,j}$ is nonzero for all $1\leq i\leq r$, $1\leq j\leq n$, then $V$ is a Mathieu-Zhao space of $B$ if and only if $\gamma_{i,1}d_1+\gamma_{i,2}d_2+\cdots+\gamma_{i,t_i}d_{t_i}\neq 0$ for some $i\in\{1,2,\ldots,r\}$ for all nonzero vectors $\tilde{d}=(d_1,d_2,\ldots,d_{t_i})$ and $d_{j_i}\in \{0,1\}$ for $1\leq j_i\leq t_i$, $t_i\in \{0,1,\ldots,n\}$ and at least one of $t_i$ is nonzero for $1\leq i\leq r$.
\end{prop}
\begin{proof}
We can assume that $\gamma_{i,j}\neq 0$ for all $1\leq j\leq t$ for some $i\in \{1,2,\ldots,r\}$ and $\gamma_{i,j}=0$ for all $1\leq i\leq r$, $t+1\leq j\leq n$ by reordering $\gamma_{i,j}$ for $1\leq i\leq r$, $1\leq j\leq n$. Then we have
$$\overbrace{0\times\cdots\times 0}^{t~columns}\times K\cdots \times K\subseteq V$$
and
$$0\times\cdots\times K\times0\times\cdots\times\overbrace{0\times\cdots\times 0}^{n-t~columns}\nsubseteq V,$$
where $t=\operatorname{max}\{t_1,t_2,\ldots,t_r\}$.

$``\Rightarrow"$ Suppose that $\gamma_{i,1}d_1+\gamma_{i,2}d_2+\cdots+\gamma_{i,t_i}d_{t_i}=0$ for some nonzero vector $\tilde{d}=(d_1,d_2,\ldots,d_{t_i})$, $d_{j_i}=0$ or 1 for $1\leq j_i\leq t_i$ for all $1\leq i\leq r$, then $e=(d_1,\ldots,d_t,0,\ldots,0)$ is an idempotent of $V$. Since $V$ is a Mathieu-Zhao space of $B$, we have that $Be=Kd_1\times\cdots\times Kd_t\times 0\times\cdots\times0\subseteq V$, which is a contradiction. Then the conclusion follows.

$``\Leftarrow"$ Let
$I=\overbrace{0\times\cdots\times 0}^{t~columns}\times K\times\cdots \times K$. Then $I$ is an ideal of $B$. We claim that $V/I$ contains no nonzero idempotent. Suppose that $e$ is a nonzero idempotent of $V/I$. Then we have $e=(e_1,e_2,\ldots,e_t)$, where $e_j=0$ or 1 for $1\leq j\leq t$. Let $\tilde{d}=(d_1,\ldots,d_t)=e\neq (0,\ldots,0)$. Then $\gamma_{i,1}d_1+\gamma_{i,2}d_2+\cdots+\gamma_{i,t_i}d_{t_i}=0$ for all $1\leq i\leq r$, which is a contradiction. It follows from Theorem 4.2 in \cite{17} that $V/I$ is a Mathieu-Zhao space of $B/I$. Then it follows from Proposition 2.7 in \cite{17} that $V$ is a Mathieu-Zhao space of $B$.
\end{proof}

\begin{rem}\label{rem3.6}
In Proposition \ref{prop3.1}, if $\gamma_{i,j}=0$ for all $1\leq i\leq r$, $1\leq j\leq n$, then $V=B$. Clearly, $V$ is a Mathieu-Zhao space of $B$.
\end{rem}

\begin{cor}\label{cor3.2}
Let $L$ and $G$ be as in Problem \ref{prob1.2} and $K$ a field of characteristic zero or a field of characteristic $p$ and $p\nmid |G|$, If $K$ is a split field for $G$ and $G$ is Abelian, then $\operatorname{Ker}L$ is a Mathieu-Zhao space of $K[G]$ if and only if $\gamma_{i,1}d_1+\gamma_{i,2}d_2+\cdots+\gamma_{i,t_i}d_{t_i}\neq 0$ for some $i\in\{1,2,\ldots,r\}$ for all nonzero vectors $\tilde{d}=(d_1,d_2,\ldots,d_{t_i})$ and $d_{j_i}\in \{0,1\}$ for $1\leq j_i\leq t_i$, $t_i\in\{0,1,\ldots,n\}$ and at least one of $t_i$ is nonzero for $1\leq i\leq r$, where $\gamma_{i,j}=\rho_j(\alpha_i)$ for all $1\leq i\leq r$, $1\leq j\leq n$ and $\rho_j$ is an irreducible representation of $G$ for $1\leq j\leq n$ and $\alpha_i$ be as in Lemma \ref{lem2.5} for $1\leq i\leq r$.
\end{cor}
\begin{proof}
Since $G$ is Abelian, we have that all the irreducible representations of $G$ are degree one. It follows from Theorem \ref{thm2.6} that $\operatorname{Ker}L\cong\{(a_1,a_2,\ldots,a_n)\in A|\sum_{j=1}^n\gamma_{i,j}a_j=0~for~all~1\leq i\leq r\}$, where $A$ is $n$ times product of $K$, $\gamma_{i,j}=\rho_j(\alpha_i)=\operatorname{Tr}\rho_j(\alpha_i)\in K$ for all $1\leq i\leq r$, $1\leq j\leq n$. Since $L\neq 0$, we have that at least one of $\gamma_{i,j}$ is nonzero for $1\leq i\leq r$, $1\leq j\leq n$. Then the conclusion follows from Proposition \ref{prop3.1}.
\end{proof}

\begin{lem}\label{lem3.3}
Let $R$ be an integral domain of characteristic $p$ and $G$ a finite Abelian group with $|G|=p^ad$, $p\nmid d$. Then every idempotent of $R[G]$ is also an idempotent of $R[\tilde{G}]$, where $G=H\times \tilde{G}$ and $|H|=p^a$. In particular, the idempotents of $R[G]$ are the same as the idempotents of $R[\tilde{G}]$.
\end{lem}
\begin{proof}
Since $G$ is a finite Abelian group, we have that $G=H\times \tilde{G}$ and $|\tilde{G}|=d$. Let $e$ be an idempotent of $R[G]$. Then $e$ can be written as
$$e=\sum_{h\in H}\alpha_hh$$
with $\alpha_h\in R[\tilde{G}]$ for each $h\in H$. Since $|H|=p^a$, we have $h^{q^m}=1$ for any $m\geq 1$, $h\in H$, where $q=p^a$. Thus, we have
$$e=e^{q^m}=\sum_{h\in H}\alpha_h^{q^m}\in R[\tilde{G}].$$
Then the conclusion follows.
\end{proof}

\begin{thm}\label{thm3.4}
Let $L$ and $G$ be as in Problem \ref{prob1.2} and $K$ a field of characteristic $p$. If $K$ is a split field for $G$ and $G$ is Abelian with $|G|=p^ad$, $p\nmid d$, then the following statements are equivalent:

$(1)$ $\operatorname{Ker}L$ is a Mathieu-Zhao space of $K[G]$.

$(2)$ $\gamma_{i,1}d_1+\gamma_{i,2}d_2+\cdots+\gamma_{i,t_i}d_{t_i}\neq 0$ for some $i\in\{1,2,\ldots,r\}$ for all nonzero vectors $\tilde{d}=(d_1,d_2,\ldots,d_{t_i})$ and $d_{j_i}\in \{0,1\}$ for $1\leq j_i\leq t_i$, $t_i\in\{0,1,\ldots,d\}$ and at least one of $t_i$ is nonzero for $1\leq i\leq r$, where $\gamma_{i,j}=\rho_j(\alpha_i)=\operatorname{Tr}\rho_j(\alpha_i)$ for $1\leq i\leq r$, $1\leq j\leq d$ and $\rho_j$ is an irreducible representation of $G/H$ for $1\leq j\leq d$, $H$ is a Sylow $p$-subgroup of $G$ and $\alpha_i$ be as in Lemma \ref{lem2.5} by replacing $G$ with $G/H$ for $1\leq i\leq r$; $l_{i,1}, l_{i,2}, \ldots, l_{i,n}$ satisfy the following equations:
\begin{equation}\label{eq3.1}
  \left\{ \begin{aligned}
  \chi_j(\tilde{g}_1^{-1})l_{i,1}+\chi_j(\tilde{g}_2^{-1})l_{i,p^a+1}+\cdots+\chi_j(\tilde{g}_d^{-1})l_{i,(d-1)p^a+1}=0\\ \chi_j(\tilde{g}_1^{-1})l_{i,2}+\chi_j(\tilde{g}_2^{-1})l_{i,p^a+2}+\cdots+\chi_j(\tilde{g}_d^{-1})l_{i,(d-1)p^a+2}=0\\
  \vdots~~~~~~~~~~~~~~~~~~~~~~~~~~~~~~~~~~~~~~~~~~~\\
  \chi_j(\tilde{g}_1^{-1})l_{i,p^a}+\chi_j(\tilde{g}_2^{-1})l_{i,2p^a}+\cdots+\chi_j(\tilde{g}_d^{-1})l_{i,dp^a}=0~~~~~~~\\
                          \end{aligned} \right.
  \end{equation}
for all $1\leq i\leq r$, $t+1\leq j\leq d$, where $\chi_j$ is the irreducible character according to $\rho_j$ for $t+1\leq j\leq d$ and $G=\cup_{k=1}^d\tilde{g}_kH$ with $\tilde{g}_1=1_{G/H}$ and $H=\{h_1,h_2,\ldots,h_{p^a}\}$ with $h_1=1_H$, $L_i(h_k)=l_{i,k}$ and $L_i(\tilde{g}_kh_q)=l_{i,(k-1)p^a+q}$ for all $1\leq i\leq r$, $1\leq k\leq d$, $1\leq q\leq p^a$ and $t=\operatorname{max}\{t_1,t_2,\ldots,t_r\}$.
\end{thm}
\begin{proof}
 Since $G$ is Abelian, we have that $G=H\times \tilde{G}$, where $\tilde{G}\cong G/H$ and $|\tilde{G}|=d$.\\ Note that
$$\gamma_{i,j}=\operatorname{Tr}\rho_j(\alpha_i)=\sum_{k=1}^d\operatorname{Tr}\rho_j(\tilde{g}_k^{-1})l_{i,(k-1)p^a+1}=\sum_{k=1}^d\chi_j(\tilde{g}_k^{-1})l_{i,(k-1)p^a+1}$$
for all $1\leq i \leq r$, $1\leq j\leq d$. Let $e_j=d^{-1}\sum_{k=1}^d\chi_j(\tilde{g}_k^{-1})\tilde{g}_k$ for $1\leq j\leq d$. Then it follows from Theorem 2.12 in \cite{4} that $e_1, e_2,\ldots,e_d$ are the primitive orthogonal idempotents of $K[\tilde{G}]$. Without loss of generality, we can assume that $\gamma_{i,j}=0$ for all $1\leq i\leq r$, $t+1\leq j\leq d$ and $\gamma_{i,j}\neq 0$ for all $1\leq j\leq t$ for some $i\in \{1,2,\ldots,r\}$ by reordering $\rho_j(\alpha_i)$ for all $1\leq i\leq r$, $1\leq j\leq d$.

$(1)\Rightarrow(2)$ It's easy to see that if $\gamma_{i,j}=0$ for all $1\leq i\leq r$, $t+1\leq j\leq d$, then $e_{t+1},\ldots,e_d$ belong to $\operatorname{Ker}(L|_{\tilde{G}})\subseteq \operatorname{Ker}L$. Thus, the ideal $I$ generated by $e_{t+1},\ldots,e_d$ belongs to $\operatorname{Ker}L$. Since $\tilde{G}$ is Abelian, it's easy to check that $e_j\tilde{g}_k=\chi_j(\tilde{g}_k)e_j$ for all $1\leq j, k\leq d$. Hence we have $e_j\tilde{g}_k\in \operatorname{Ker}L$ for all $t+1\leq j\leq d$, $1\leq k\leq d$. Note that $e_jh_q\in \operatorname{Ker}L$ for all $t+1\leq j\leq d$, $1\leq q\leq p^a$. Then we have equations \eqref{eq3.1} for all $1\leq i\leq r$, $t+1\leq j\leq d$. It follows from Proposition \ref{prop2.3} that $\operatorname{Ker}(L|_{\tilde{G}})$ is a Mathieu-Zhao space of $K[\tilde{G}]$. That is, $\operatorname{Ker}(L|_{G/H})$ is a Mathieu-Zhao space of $K[G/H]$. Since $p\nmid |G/H|$, the conclusion follows from Corollary \ref{cor3.2}.

$(2)\Rightarrow(1)$ If $\gamma_{i,j}=0$ for all $1\leq i\leq r$, $t+1\leq j\leq d$, then $e_{t+1},\ldots,e_d\in \operatorname{Ker}(L|_{\tilde{G}})\subseteq \operatorname{Ker}L$. It's easy to check that $e_j\tilde{g}_k=\chi_j(\tilde{g}_k)e_j$ and $e_j\tilde{g}_kh_q=\chi_j(\tilde{g}_k)e_jh_q$ for all $t+1\leq j\leq d$, $1\leq k\leq d$, $1\leq q\leq p^a$. Therefore, we have $I\subseteq \operatorname{Ker}L$, where $I$ is an ideal generated by $e_{t+1},\ldots,e_d$. Since $e_1,\ldots,e_d$ are the primitive orthogonal idempotents of $K[\tilde{G}]$ and there are $2^d$ idempotents in $K[\tilde{G}]$, we have that any idempotent of $K[\tilde{G}]$ is a sum of some of the $e_j$ for $1\leq j\leq d$. Note that the condition that $\gamma_{i,1}d_1+\gamma_{i,2}d_2+\cdots+\gamma_{i,t_i}d_{t_i}\neq 0$ for some $i\in\{1,2,\ldots,r\}$ for all nonzero vectors $\tilde{d}=(d_1,d_2,\ldots,d_{t_i})$ and $d_{j_i}\in \{0,1\}$ is equivalent to that any sum of some of the $e_j$ is not in $\operatorname{Ker}(L|_{\tilde{G}})$ except zero for all $1\leq j\leq t$ . Hence any sum of some of the $e_j$ is not in $\operatorname{Ker}(L|_{\tilde{G}})$ for all $1\leq j\leq d$ if the sum contains $e_{j_0}$ for some $j_0\in \{1,2,\ldots,t\}$. Thus, any sum of some of the $e_j$ is not in $\operatorname{Ker}L$ for all $1\leq j\leq d$ if the sum contains $e_{j_0}$ for some $j_0\in \{1,2,\ldots,t\}$. Otherwise, the sum of $e_j$ belong to $\operatorname{Ker}L\cap K[\tilde{G}]=\operatorname{Ker}(L|_{\tilde{G}})$ for $1\leq j\leq d$, which is a contradiction. It follows from Lemma \ref{lem3.3} that $K[G]$ and $K[\tilde{G}]$ have the same idempotents. Hence $\operatorname{Ker}L/I$ has no nonzero idempotent. It follows from Theorem 4.2 in \cite{17} that $\operatorname{Ker}L/I$ is a Mathieu-Zhao space of $K[G]/I$. Hence it follows from Proposition 2.7 in \cite{17} that $\operatorname{Ker}L$ is a Mathieu-Zhao space of $K[G]$.
\end{proof}

\begin{rem}\label{rem3.5}
If $G$ is cyclic in Theorem \ref{thm3.4}, then all the primitive orthogonal idempotents of $K[G]$ are $e_j=d^{-1}(1+(\xi^{d-1})^{j-1}\tilde{g}+\cdots+\xi^{j-1}\tilde{g}^{d-1})$ for $1\leq j\leq d$, where $\xi$ is a $d$-th root of unity and $\tilde{G}$ is generated by $\tilde{g}$, where $\tilde{G}$ be as in Theorem \ref{thm3.4}.
\end{rem}

{\bf{Acknowledgement}}:  The author is very grateful to professor Wenhua Zhao for some useful suggestions. She is also grateful to the Department of Mathematics of Illinois State University, where this paper was partially finished, for hospitality during her stay as a visiting scholar.

\end{document}